\documentclass[final,1p,times]{elsarticle}

\usepackage{amssymb}
 \usepackage{amsthm}
\usepackage{amscd}
\usepackage{amsmath}
\usepackage{amsfonts}
\usepackage{amssymb}
\usepackage{graphicx}

\numberwithin{equation}{section}

\newcommand{\f}{\frac}

\newcommand{\be}{\begin{equation}}

\newcommand{\ee}{\end{equation}}
\newcommand{\bea}{\begin{eqnarray}}
\newcommand{\eea}{\end{eqnarray}}
\newcommand{\bna}{\begin{eqnarray*}}
\newcommand{\ena}{\end{eqnarray*}}

\renewcommand{\le}{\left}
\newcommand{\ri}{\right}

\journal{***}

\begin{document}

\begin{frontmatter}

\title{A new proof of subcritical Trudinger-Moser inequalities on the whole Euclidean space}

\author{Yunyan Yang}
 \ead{yunyanyang@ruc.edu.cn}
 \author{Xiaobao Zhu}
 \ead{zhuxiaobao@ruc.edu.cn}
\address{ Department of Mathematics,
Renmin University of China, Beijing 100872, P. R. China}

\begin{abstract}
In this note, we give a new proof of subcritical Trudinger-Moser
inequality on $\mathbb{R}^n$. All the existing proofs on this
inequality are based on the rearrangement argument with respect to
functions in the Sobolev space $W^{1,n}(\mathbb{R}^n)$. Our method
avoids this technique and thus can be used in the Riemannian manifold
case and in the entire Heisenberg group.
\end{abstract}


\end{frontmatter}

\section{Introduction}

 It was proved by Cao \cite{Cao}, Panda \cite{Panda} and do
 \'O \cite{doo1} that\\

 \noindent {\bf Theorem A} {\it
Let $\alpha_{n}=n\omega_{n-1}^{\frac{1}{n-1}}$, where $\omega_{n-1}$ is the measure of the unit sphere in $\mathbb{R}^{n}$.
Then for any $\alpha<\alpha_{n}$ there holds
\begin{align}\label{e1-1}
  \sup_{u\in W^{1,n}(\mathbb{R}^n),\,\,\int_{\mathbb{R}^{n}}(|\nabla \
    u|^{n}+|u|^{n})dx\leq1}\int_{\mathbb{R}^{n}}\le(e^{\alpha|u|^{\frac{n}{n-1}}}-\sum_{k=0}^{n-2}\frac{\alpha^{k}|u|^{\frac{nk}{n-1}}}{k!}\ri)dx<\infty.
\end{align}
}

This result has various extensions, among which we mention Adachi and Tanaka \cite{AT}, Ruf \cite{Ruf}, Li-Ruf \cite{LiRuf},
Adimurthi-Yang \cite{Adi-Yang}. To the authors' knowledge, all the existing proofs of such an inequality are based on rearrangement argument with respect to
functions in the Sobolev space $W^{1,n}(\mathbb{R}^n)$. The purpose of this short note is to provide a new method to reprove Theorem A. Namely, we use a technique
of the analogy of unity decomposition. More precisely, for any $u\in W^{1,n}(\mathbb{R}^n)$, we first take a cut-off function
$\phi_i\in C_0^\infty(B_R(x_i))$ such that $0\leq\phi_i\leq 1$ on $B_R(x_i)$, $\phi_i\equiv 1$ on $B_{R/2}(x_i)$.
Then, using the usual Trudinger-Moser inequality \cite{Moser,Pohozaev,Trudinger} for bounded domain, we prove a key estimate
\be\label{key}\int_{\mathbb{R}^{n}}\le(e^{\alpha|\phi_iu|^{\frac{n}{n-1}}}-\sum_{k=0}^{n-2}\frac{\alpha^{k}|\phi_iu|
    ^{\frac{nk}{n-1}}}{k!}\ri)dx\leq \int_{\mathbb{R}^n}|\nabla (\phi_iu)|^ndx\ee
under the condition that
$$\int_{\mathbb{R}^{n}}|\nabla (\phi_i u)|^{n}dx\leq1.$$
The power of (\ref{key}) is evident. It permits us to approximate $u$ by $\sum_i\phi_i u$, where every $\phi_i$ is supported in $B_{R}(x_i)$,
$\mathbb{R}^n=\cup_{i=1}^\infty B_{R/2}(x_i)$, and any fixed $x\in\mathbb{R}^n$ belongs to at most $c(n)$ balls $B_R(x_i)$ for some universal constant $c(n)$. 
 If we further take $\phi_i$ such that $|\nabla\phi_i|\leq 4/R$. 
Note that for any $\epsilon>0$ there exists a constant $C(\epsilon)$ such that
$$\int_{\mathbb{R}^{n}}|\nabla (\phi_i u)|^{n}dx\leq (1+\epsilon)\int_{\mathbb{R}^{n}}|\nabla  u|^{n}dx+\f{C(\epsilon)}{R^n}
\int_{\mathbb{R}^{n}}|u|^{n}dx.$$
Selecting $\epsilon>0$ sufficiently small and $R>0$ sufficiently large, we get the desired result. \\

Similar idea was used by the first named author to deal with similar problems on complete Riemannian manifolds \cite{manifold}
or the entire Heisenberg group \cite{Y-Heisenberg}. Note that due to the complicated geometric structure, we have not obtained Theorem A on manifolds,
but a weaker result. Namely\\

\noindent{\bf Theorem B} {\it  Let $(M,g)$ be a complete noncompact Riemannian $n$-manifold.
 Suppose that its Ricci curvature has lower bound, namely
 ${\rm Rc}_{(M,g)}\geq Kg$ for some constant $K\in\mathbb{R}$, and its injectivity radius is strictly positive,
 namely ${\rm inj}_{(M,g)}\geq i_0$ for some constant $i_0>0$. Then
 we have

 \noindent $(i)$ for any $0\leq\alpha<\alpha_n$ there
  exists positive constants $\tau$ and $\beta$ depending only on $n$, $\alpha$, $K$ and $i_0$ such that
  \be\label{Tm-hold}\sup_{u\in W^{1,n}(M),\,\|u\|_{1,\tau}\leq 1}
 \int_M\le(e^{\alpha|u|^{\f{n}{n-1}}}-\sum_{k=0}^{n-2}\f{\alpha^k|u|^{\f{nk}{n-1}}}{k!}\ri)dv_g
 \leq \beta,\ee
 where \be\label{1tau}\|u\|_{1,\tau}=\le(\int_M|\nabla_g u|^ndv_g\ri)^{1/n}+
 \tau\le(\int_M|u|^ndv_g\ri)^{1/n}.\ee
 As a consequence, $W^{1,n}(M)$ is embedded in $L^q(M)$ continuously
 for all $q\geq n$;

 \noindent $(ii)$ for any $\alpha>\alpha_n$ and any $\tau>0$, the supremum in (\ref{Tm-hold}) is
 infinite;

 \noindent $(iii)$ for any $u\in W^{1,n}(M)$ and any $\alpha>0$, the integrals in (\ref{Tm-hold}) are still finite.\\}

 We say more words about this method. For Sobolev inequalities on complete noncompact Riemannian manifolds, unity decomposition was employed by
 Hebey et al. \cite{Hebey}. In the case of Trudinger-Moser inequality, it is not evidently applicable. We are lucky to find its analogy 
 (\cite{manifold}, Lemma 4.1).

\section{Preliminary lemmas}

We first give a local estimate concerning the Trudinger-Moser functional. Precisely we have\\

\noindent{\bf Lemma 1} {\it 
For any $x_0\in \mathbb{R}^n$ and any $u\in W_{0}^{1,n}(B_{R}(x_0))$, $\int_{B_{R}(x_0)}|\nabla u|^{n}dx\leq1$, we have
\begin{align}\label{e2-1}
  \int_{B_{R}(x_0)}\le(e^{\alpha_{n}|u|^{\frac{n}{n-1}}}
                -\sum_{k=0}^{n-2}\frac{\alpha_{n}^{k}|u|^{\frac{nk}{n-1}}}{k!}\ri)dx
                \leq C(n)R^{n}\int_{B_{R}(x_0)}|\nabla u|^{n}dx,
\end{align}
where $C(n)$ is a constant depending only on $n$.}\\

\noindent{\it Proof.} Essentially this is the same as (\cite{manifold}, Lemma 4.1). For reader's convenience we give the
details here.  It is well known \cite{Moser,Pohozaev,Trudinger} that
\begin{align}\label{e2-2}
  \sup_{u\in W_{0}^{1,n}(B_{R}(x_0)),
  \int_{B_{R}(x_0)}|\nabla u|^{n}dx\leq1}\int_{B_{R}(x_0)}e^{\alpha_{n}|u|^{\frac{n}{n-1}}}
                dx
                \leq C(n)R^{n}.
\end{align}
Letting $\widetilde{u}=\frac{u}{||\nabla u||_{L^{n}({B_{R}(x_0)})}}$ for any $u\in W_{0}^{1,n}(B_{R}(x_0))\setminus\{0\}$, we have
\bea\label{e2-3}
    \int_{B_{R}(x_0)}\le(e^{\alpha_{n}|\widetilde{u}|^{\frac{n}{n-1}}}
                -\sum_{k=0}^{n-2}\frac{\alpha_{n}^{k}|\widetilde{u}|^{\frac{nk}{n-1}}}{k!}\ri)dx
&\geq&\frac{1}{||\nabla u||_{L^{n}({B_{R}(x_0)})}}\int_{B_{R}(x_0)}\sum_{k=n-1}^{\infty}\frac{\alpha_{n}^{k}|u|^{\frac{nk}{n-1}}}{k!}dx\nonumber\\
   &=&\frac{1}{||\nabla u||_{L^{n}({B_{R}(x_0)})}}\int_{B_{R}(x_0)}\le(e^{\alpha_{n}|u|^{\frac{n}{n-1}}}
                -\sum_{k=0}^{n-2}\frac{\alpha_{n}^{k}|u|^{\frac{nk}{n-1}}}{k!}\ri)dx.\quad
\eea
Combining (\ref{e2-2}) and (\ref{e2-3}), we get the desired result. $\hfill\Box$\\

Also we need a covering lemma of $\mathbb{R}^{n}$, see for example (\cite{Hebey}, Lemma 1.6). \\

\noindent{\bf Lemma 2} {\it 
For any $R>0$, there exists a sequence $\{x_{i}\}_{i=1}^{\infty}\subset\mathbb{R}^{n}$
such that\\
$(i)$ $\cup_{i=1}^{\infty}B_{{R}/{2}}(x_{i})=\mathbb{R}^{n}$;\\
$(ii)$ $\forall i\neq j,\, B_{{R}/{4}}(x_{i})\cap B_{{R}/{4}}(x_{j})=\varnothing$;\\
$(iii)$ $\forall x\in\mathbb{R}^{n}$, $x$ belongs to at most $N$ balls $B_{R}(x_{i})$ for some integer $N$.}

\section{Proof of Theorem A}

We shall obtain a global inequality (\ref{e1-1}) by gluing local estimates (\ref{e2-1}). \\

\noindent{\it Proof of Theorem A.} Let $R>0$ to be determined later.
Let $\phi_{i}$ be the cut-off function satisfies the following conditions:
$(i)$ $\phi_{i}\in C_{0}^{\infty}(B_{R}(x_{i}))$;
   $(ii)$ $0\leq\phi_{i}\leq1$ on $B_R(x_i)$ and $\phi_{i}\equiv 1$ on $B_{{R}/{2}}(x_{i})$;
   $(iii)$ $|\nabla \phi_{i}(x)|\leq{4}/{R}$.
For $u\in W^{1,n}(\mathbb{R}^{n})$ satisfying
\begin{align}\label{e3-1}
\int_{\mathbb{R}^{n}}(|\nabla u|^{n}+|u|^{n})dx\leq1,
\end{align}
we have $\phi_{i}u\in W_{0}^{1,n}(B_{R}(x_{i}))$, using Cauchy inequality with $\epsilon$ term we obtain
\begin{align}\label{e3-2}
\int_{B_{R}(x_{i})}|\nabla(\phi_{i}u)|^{n}dx\leq&(1+\epsilon)\int_{B_{R}(x_{i})}\phi_{i}^{n}|\nabla u|^{n}dx
         +C(\epsilon)\int_{B_{R}(x_{i})}|\nabla \phi_{i}|^{n}|u|^{n}dx\nonumber\\
         \leq&(1+\epsilon)\int_{B_{R}(x_{i})}|\nabla u|^{n}dx
         +\frac{C(\epsilon)}{R^{n}}\int_{B_{R}(x_{i})}|u|^{n}dx\nonumber\\
         \leq&(1+\epsilon)\int_{B_{R}(x_{i})}(|\nabla u|^{n}+|u|^{n})dx,
\end{align}
where in the last inequality we choose a sufficiently large $R$ to make sure $\frac{C(\epsilon)}{R^{n}}\leq(1+\epsilon)$.
Let $\alpha_{\epsilon}=\frac{\alpha_{n}}{(1+\epsilon)^{1/(n-1)}}$ and $\widetilde{\phi_{i}u}=\frac{\phi_{i}u}{(1+\epsilon)^{1/n}}$.
Noting that $\widetilde{\phi_{i}u}\in W_{0}^{1,n}(B_{R}(x_{i}))$, we have by (\ref{e3-2}) and Lemma 1 
\begin{align}\label{e3-3}
    \int_{B_{\frac{R}{2}}(x_{i})}\le(e^{\alpha_{\epsilon}|u|^{\frac{n}{n-1}}}-\sum_{k=0}^{n-2}\frac{\alpha_{\epsilon}^{k}|u|^{\frac{nk}{n-1}}}{k!}\ri)dx
\leq&\int_{B_{R}(x_{i})}\le(e^{\alpha_{\epsilon}|\phi_{i}u|^{\frac{n}{n-1}}}-\sum_{k=0}^{n-2}\frac{\alpha_{\epsilon}^{k}|\phi_{i}u|^{\frac{nk}{n-1}}}{k!}\ri)dx\nonumber\\
   =&\int_{B_{R}(x_{i})}\le(e^{\alpha_{n}|\widetilde{\phi_{i}u}|^{\frac{n}{n-1}}}-\sum_{k=0}^{n-2}\frac{\alpha_{n}^{k}|\widetilde{\phi_{i}u}|^{\frac{nk}{n-1}}}{k!}\ri)dx
   \nonumber\\
\leq&C(n)R^{n}\int_{B_{R}(x_{i})}|\nabla(\widetilde{\phi_{i}u})|^{n}dx\nonumber\\
\leq&C(n)R^{n}\int_{B_{R}(x_{i})}(|\nabla u|^{n}+|u|^{n})dx.
\end{align}
By Lemma 2 and (\ref{e3-3}), we have
\begin{align}\label{e3-4}
   \int_{\mathbb{R}^{n}}\le(e^{\alpha_{\epsilon}|u|^{\frac{n}{n-1}}}-\sum_{k=0}^{n-2}\frac{\alpha_{\epsilon}^{k}|u|^{\frac{nk}{n-1}}}{k!}\ri)dx
\leq&\int_{\cup_{i=1}^{\infty}B_{\frac{R}{2}}(x_{i})}\le(e^{\alpha_{\epsilon}|u|^{\frac{n}{n-1}}}-\sum_{k=0}^{n-2}\frac{\alpha_{\epsilon}^{k}|u|^{\frac{nk}{n-1}}}{k!}\ri)dx\nonumber\\
\leq&\sum_{i=1}^{\infty}\int_{B_{\frac{R}{2}}(x_{i})}\le(e^{\alpha_{\epsilon}|u|^{\frac{n}{n-1}}}-\sum_{k=0}^{n-2}\frac{\alpha_{\epsilon}^{k}|u|^{\frac{nk}{n-1}}}{k!}\ri)dx\nonumber\\
\leq&\sum_{i=1}^{\infty}C(n)R^{n}\int_{B_{R}(x_{i})}(|\nabla u|^{n}+|u|^{n})dx\nonumber\\
\leq&C(n)R^{n}N\int_{\mathbb{R}^{n}}(|\nabla u|^{n}+|u|^{n})dx\nonumber\\
\leq&C(n)R^{n}N.
\end{align}
For any $\alpha<\alpha_n$, we can choose $\epsilon>0$ sufficiently small such that $\alpha<\alpha_\epsilon$. 
This ends the proof of Theorem A. $\hfill\Box$

\section{Concluding remarks}
Using the same idea to prove Theorem A, we can also prove the subcritical Adams inequality in $\mathbb{R}^n$ \cite{Adams,Ruf-Sani,Yang-JDE},
which strengthen (\cite{manifold}, Theorem 2.6). Since
the proof is completely analogous to our proof of Theorem A, we leave it to the reader.\\

{\bf Acknowledgement.} This work is supported by the NSFC 11171347.


\begin{thebibliography}{00}


\bibitem{AT} S. Adachi, K. Tanaka, Trudinger type inequalities in $\mathbb{R}^N$
and their best exponents, Proc. Amer. Math. Soc. 128 (2000)
2051-2057.

\bibitem{Adams} D. Adams, A sharp inequality of J. Moser for
higher order derivatives, Ann. Math. 128 (1988) 385-398.


\bibitem{Adi-Yang} Adimurthi, Y. Yang, An interpolation of Hardy inequality
and Trudinger-Moser inequality in $\mathbb{R}^N$ and its
applications, Internat. Mathematics Research Notices 13 (2010)
2394-2426.



\bibitem{Cao} D. Cao, Nontrivial solution of semilinear elliptic
equations with critical exponent in $\mathbb{R}^2$, Commun. Partial
Differential Equations 17 (1992) 407-435.




\bibitem{doo1} J. M. do \'O, $N$-Laplacian equations in
$\mathbb{R}^N$ with critical growth, Abstr. Appl. Anal. 2 (1997)
301-315.






\bibitem{Hebey} E. Hebey, Sobolev spaces on Riemannian maifolds,
Lecture notes in mathematics 1635, Springer, 1996.


\bibitem{LiRuf} Y. Li, B. Ruf, A sharp Trudinger-Moser type
inequality for unbounded domains in $\mathbb{R}^N$, Ind. Univ. Math.
J. 57 (2008) 451-480.


\bibitem{Moser}  J. Moser, A sharp form of an inequality by
N.Trudinger, Ind. Univ. Math. J. 20 (1971) 1077-1091.

\bibitem{Panda} R. Panda, Nontrivial solution of a quasilinear
elliptic equation with critical growth in $\mathbb{R}^n$, Proc.
Indian Acad. Sci. (Math. Sci.)  105 (1995) 425-444.

\bibitem{Pohozaev} S. Pohozaev, The Sobolev embedding in the special case
$pl=n$, Proceedings of the technical scientific conference on
advances of scientific reseach 1964-1965, Mathematics sections,
158-170, Moscov. Energet. Inst., Moscow, 1965.


\bibitem{Ruf} B. Ruf, A sharp Trudinger-Moser type inequality for
unbounded domains in $\mathbb{R}^2$, J. Funct. Anal. 219 (2005)
340-367.

\bibitem{Ruf-Sani} B. Ruf, F. Sani, Sharp Adams-type inequalities in
$\mathbb{R}^n$. Trans. Amer. Math. Soc. (In press).


\bibitem{Trudinger} N. S. Trudinger, On embeddings into Orlicz spaces and
some applications, J. Math. Mech. 17 (1967) 473-484.


\bibitem{manifold} Y. Yang, Trudinger-Moser inequalities on complete noncompact Riemannian manifolds,
J. Funct. Anal. 263 (2012) 1894-1938.



\bibitem{Yang-JDE} Y. Yang, Adams type inequalities and related elliptic partial differential equations in dimension four,
J. Differ. Equations 252 (2012) 2266-2295.

\bibitem{Y-Heisenberg} Y. Yang, Trudinger-Moser inequalities on the entire Heisenberg group, arXiv:1201.2993.


\end{thebibliography}
\end{document}